\documentclass[10pt,fullpage]{article}
\usepackage{amsmath}
\usepackage{amsfonts,amssymb}
\usepackage{amsthm}
\usepackage{graphics,graphicx}
\usepackage{hyperref}
\usepackage[right = 2.5cm, left=2.5cm, top = 2.5cm, bottom =2.5cm]{geometry}
\newcommand{\Real}{\mathbb R}
\newcommand{\norm}[1]{\|#1\|}
\newcommand{\abs}[1]{\left\vert#1\right\vert}
\newcommand{\set}[1]{\left\{#1\right\}}

\newcommand{\grad}{\nabla}

\newcommand{\K}{\mathcal{K}}

\newcommand{\wknorm}[2]{\norm{#1}_{L^{#2,\infty}}}

\newcommand{\F}{\mathcal{F}}

\newtheorem{theorem}{Theorem}
\theoremstyle{remark}
\newtheorem{remark}{Remark}

\theoremstyle{theorem}
\newtheorem{prop}{Proposition}

\theoremstyle{definition}

\theoremstyle{lemma}
\newtheorem{lemma}{Lemma}
\theoremstyle{corollary}

\begin{document}

\title{Global Minimizers for Free Energies of Subcritical Aggregation Equations with Degenerate Diffusion} 

\author{Jacob Bedrossian \footnote{\textit{jacob.bedrossian@math.ucla.edu}, University of California-Los Angeles, Department of Mathematics}}

\date{}

\maketitle

\begin{abstract}
We prove the existence of non-trivial global minimizers of a class of free energies related to aggregation equations with degenerate diffusion on $\Real^d$.
Such equations arise in mathematical biology as models for organism group dynamics which account for competition between the tendency to aggregate into groups and nonlinear diffusion to avoid over-crowding. The existence of non-zero optimal free energy stationary solutions representing coherent groups in $\Real^d$ is therefore of interest. 
The primary contribution is the investigation of a notion of criticality associated with the global minimizer problem.
The notion arises from the scaling of diffusive and aggregative forces as mass spreads and is shown to dictate the existence, and sometimes non-existence, of global minimizers. 
\end{abstract} 

\section{Introduction} 
We consider global minimizers of the non-convex free energy
\begin{align}
\F(u) = \int_{\Real^d} \Phi(u) dx - \frac{1}{2}\int\int_{\Real^d\times\Real^d} u(x)u(y) \K(x-y) dxdy := S(u) - \frac{1}{2}\mathcal{W}(u) \label{def:F},
\end{align}
with $u \in L_+^1(\Real^d) := \set{u \in L^1(\Real^d): u\geq 0}$ such that $\norm{u}_1 := M$ over $\Real^d$ with $d \geq 2$. We are interested in determining for which choices of $\K$, $\Phi$ and $M$ there exist global minimizers with mass $M$. We restrict to the case where $\Phi$ and $\K$ are non-negative. We refer to $S(u)$ as the entropy and $\mathcal{W}(u)$ as the interaction energy.  Free energies such as \eqref{def:F} arise in the study of aggregation equations and Patlak-Keller-Segel models with degenerate diffusion
\cite{BRB10,TopazBertozziLewis06,Burger07,LiRodrigo09,Blanchet09,Slepcev08,BertozziSlepcev10,SugiyamaNA05,SugiyamaADE07,SugiyamaDIE06,SugiyamaDIE07}. A typical example is 
\begin{equation}
u_t + \grad \cdot (u \grad \K \ast u) = \Delta u^m \label{def:ADD}
\end{equation} 
for $m > 1$ and $u(t) \in L_+^1(\Real^d)$.  These equations are formally a gradient flow for \eqref{def:F} with $\Phi(u) = \frac{1}{m-1}u^m$ under the Euclidean Wasserstein distance \cite{Otto01,AmbrosioGigliSavare,BlanchetCalvezCarrillo08}, however as \eqref{def:F} is not displacement convex in the sense of \cite{McCann97}, the established theory does not apply. 
 Numerical simulations indicate that for certain choices of $\K$ and $m$, there are compactly supported stationary solutions to \eqref{def:ADD} which also appear to be attractors \cite{TopazBertozzi04,TopazBertozziLewis06}. 
The purpose of this work is to provide a rigorous analysis of the existence of these stationary solutions. 
Specifically, we give sufficient conditions for when there exist non-trivial global minimizers to \eqref{def:F}, which
 are formally stationary solutions to \eqref{def:ADD} with optimal free energy. We only consider the case $\K$ radially symmetric non-increasing, which in particular implies that the non-local interaction is purely attractive. 
The minimizer problem without diffusion and more general $\K$ has been considered elsewhere and is generally of a different flavor (see
 e.g. \cite{FellnerRaoul1,FellnerRaoul2}). 

A typical approach to proving the existence of a global minimizer is to take a suitably strong limit of infimizing sequences. 
On bounded domains this may be carried through without issue, however on $\Real^d$ one must deal with the possibility that the sequence can lose some or all of the mass in the limit. 
To overcome this, we use concentration compactness arguments to establish tightness of infimizing sequences. Indeed, \eqref{def:F} was considered
in the original work of Lions \cite{LionsCC84}. 
In \cite{LionsCC84}, the existence of minimizers was proved without the use of symmetry arguments in a number of cases using the concept of strict sub-additivity. 
Lions shows that this is effective for kernels with slow decay at infinity, but cases such as $\K \in L^1$ are not treated in detail.
Proposition \ref{prop:Lions} below contains the relevant results.
Therefore, the primary contribution of this work is that it applies to kernels with fast decay. 
In contrast to \cite{LionsCC84}, we rely on symmetrization arguments and a scaling analysis described below, as we could not verify strict sub-additivity in these cases. 

We consider $\K$ non-negative, radially symmetric non-increasing such that $\K \in L_{loc}^{p,\infty}(\Real^d) \cap L^{\hat{p}}(\Real^d \setminus B_1(0))$ for some $1 < p < \infty$ and $1 \leq \hat{p} < \infty$. We define $m^\star = (p+1)/p$. 
Note that we are only considering the case $m^\star > 1$. The methods below allow us to define $m^\star = 1$ if $\K \in L^\infty$, however, if $\K$ has (for instance) a
logarithmic singularity, the methods require us to define $m^\star > 1$, in contrast to \cite{BRB10}.  
Recall the following well known inequality which can be considered a generalization of the Hardy-Littlewood-Sobolev inequality: for all $f \in L^1\cap L^{m^\star}$ and $\delta > 0$, 
\begin{equation}
\int\int f(x) f(y) \K(x-y) \mathbf{1}_{B_\delta}(\abs{x-y}) dx dy \leq C_0\wknorm{\K\mathbf{1}_{B_\delta}}{p} \norm{f}_{\frac{2p}{2p-1}}^2 \leq  C_0\wknorm{\K\mathbf{1}_{B_\delta}}{p} \norm{f}_1^{2-m^\star}\norm{f}_{m^\star}^{m^\star}, \label{ineq:GHLS}
\end{equation}
where $C_0$ is a constant that depends on $p$ and $d$. Similar inequalities are used in \cite{Blanchet09,LionsCC84,BRB10,Lieb83}.
Let $\mathcal{Y}_M := \set{u \in L_+^1 \cap L^{m^\star}: \norm{u}_1 = M}$ and denote $\inf_{u \in \mathcal{Y}_M}\F(u) := I_M$. 
Note that by considering a suitable sequence weak$^\star$ converging to zero, necessarily
$I_M \leq 0$ for all $M \geq 0$.  We take $\Phi(u)$ to be strictly convex, non-negative and we assume $\K,\Phi$ and $M$ satisfy,
\begin{equation}
\lim\inf_{z \rightarrow \infty} \frac{\Phi(z)}{z^{m^\star}} > \frac{1}{2}C_0\wknorm{\K \mathbf{1}_{B_\delta}}{p}M^{2-m^\star}, \label{cond:critmass}
\end{equation}
for some $\delta > 0$. Note that the quantity on the left hand side of course need not be finite.  
Condition \eqref{cond:critmass} implies that the problem is not supercritical in the sense of \cite{BRB10} and if it is critical, then the mass is assumed to be strictly less than the critical mass.
The purpose of this assumption is to ensure $- \infty < I_M$ and that sequences $\set{u_n} \subset \mathcal{Y}_M$ with $\sup \F(u_n) < \infty$ 
have uniformly bounded entropy and $L^{m^\star}$ norms (see for instance \cite{LionsCC84} or \cite{BRB10}). 
However, we point out that kernels considered here are allowed to be more singular and have less decay than those 
considered in \cite{BRB10}. Moreover, we do not need regularity assumptions on $\K$. 
In \cite{BRB10}, the primary purpose of the limitations was to ensure solutions to \eqref{def:ADD} were unique in $L_+^1 \cap L^\infty$. 
As in \cite{LionsCC84}, we also require that $\lim_{z \rightarrow 0} \Phi(z)z^{-1} = 0$.

Before stating the main theorem, we describe the notion of criticality associated with \eqref{def:F} which dictates the main existence results in this work.
The profile decomposition in \cite{LionsCC84} applied to sequences in $\mathcal{Y}_M$ with bounded $L^{m^\star}$ norm
suggests that the primary difficulties for proving $I_M$ is attained for some $u^\star \in \mathcal{Y}_M$ will be ensuring the mass of infinimizing sequences does not 
split apart or vanish. 
Naturally, this leads to considering the relative balance of the entropy and interaction energy 
as the mass of a sequence in $\mathcal{Y}_M$ spreads out. 
For simplicity, consider the case $S(u) = \frac{1}{m-1}\int u^m dx$ with $m > m^\star \geq 1$, which arises from 
power-law degenerate diffusion, and suppose $\K \in L^1$. Let $u(x) \in L^1_+\cap L^{m}$ and define $u_\lambda(x) = \lambda^{d}u(\lambda x)$. Then, 
\begin{equation*}
\F(u_\lambda) = \frac{\lambda^{dm - d}}{m-1}\norm{u}_m^m - \frac{\lambda^{d}}{2}\int\int u(x) u(y)\frac{1}{\lambda^{d}}\K\left(\frac{x-y}{\lambda}\right) dx dy. 
\end{equation*}
The key observation is that as $\lambda \rightarrow 0$, the second term behaves like $(\lambda^d/2)\norm{\K}_1\norm{u}_2^2$.
If $m > 2$ then for sufficiently small $\lambda$, $\F(u_\lambda) < 0$, whereas if $m < 2$, this is no longer true. 
The case $m = 2$ is in some sense critical, since the entropy and interaction energy scale the same as $\lambda \rightarrow 0$, and the sign in the limit only depends on the value of $\norm{\K}_1$. 
This scaling analysis
is an important step (Lemma \ref{lem:infF}) to the proof of the main theorem below. 
We note that $m > 2$ are the exponents for which the equation \eqref{def:ADD} may be formally
re-written as a regularized non-local interface problem \cite{Slepcev08}.

\begin{theorem} \label{thm:Minimizers}
Let $m^\star > 1$, $\K$ non-trivial, $M > 0$ and the hypotheses described above be satisfied. 
Suppose $\exists\, \chi, \;  0 \leq \chi < \infty$ such that
\begin{equation}
\Phi(z) = \chi z^2 + o(z^2),\; z \rightarrow 0, \label{def:chi}
\end{equation} 
and additionally suppose that either of the following holds
\begin{itemize}
\item[(i)] $\chi = 0$. 
\item[(ii)] $0 < \chi < \infty$ and $2\chi < \norm{\K}_1 \leq \infty$.
\end{itemize} 
Then $I_M < 0$ and there exists a radially symmetric non-increasing $u^\star \in \mathcal{Y}_M$ such that $\F(u^\star) = I_M$. 
\end{theorem} 

Lions in \cite{LionsCC84} states results concerning the case when \eqref{def:chi} does not hold for any $\chi$ (see \textit{(iii)} below).
The contrast between the results shows how the scaling analysis manifests in the minimizer problem.
As mentioned above, Lions also presents results which apply to kernels with slower decay at infinity (see \textit{(iv)} below). 
Note that when kernels have slow decay, it is possible to prove the existence of minimizers in a wider variety of cases than Theorem \ref{thm:Minimizers} provides.
The following statement is a consequence of Corollary II.1 and Theorem II.1 in \cite{LionsCC84}: 
\begin{prop}(Lions \cite{LionsCC84}) \label{prop:Lions}
Suppose either of the following holds
\begin{itemize}
\item[(iii)] The mass $M$ is sufficiently large and there exists $1 < \nu < 2$ such that for all $t \geq 1$ and $z > 0$, 
\begin{equation*}
\Phi(tz) \leq t^\nu \Phi(z).
\end{equation*}
\item[(iv)] $\exists \, \alpha \in (0,d)$ such that
$\forall t \geq 1$ and $\forall \xi \in \Real^d$,
\begin{equation}
\K(t\xi) \geq t^{-\alpha}\K(\xi) \label{ineq:decay_condition}
\end{equation}
and 
\begin{equation*}
\lim_{z \rightarrow 0}\frac{\Phi(z)}{z^{1 + \alpha/d}} = 0. 
\end{equation*}
\end{itemize}
Then $I_M < 0$ and there exists a radially symmetric non-increasing $u^\star \in \mathcal{Y}_M$ such that $\F(u^\star) = I_M$.
\end{prop}

\begin{remark}
Due to the lack of convexity,
to the author's knowledge, uniqueness is largely unresolved except when $\K$ is the Newtonian potential or has similar special properties \cite{YaoKim10,LiebYau87}. In critical cases, it is known to be not unique at the critical mass \cite{Blanchet09}. 
\end{remark}

\begin{remark}
In applications, $\Phi$ is often negative near zero \cite{CarrilloEntDiss01}. However, for problems with degenerate diffusion, one can show $S(u) \gtrsim -M$, and the methods here apply simply by modifying \eqref{def:F} by an irrelevant constant depending on $M$. However, for problems with more general $\Phi$, such as the Boltzmann entropy $\Phi(u) = u\log u$, the methods here would have to be
modified.
Similarly, modifications would have to be made to treat potentials $\K$ which are unbounded from below, such as the logarithmic potential.  
\end{remark}

We now point out that the condition $\norm{\K}_1 > 2\chi$ in (ii) is close to sharp.
\begin{prop}
Let $S(u) = \norm{u}_2^2$ and $\norm{\K}_1 < 2$. Then for all $M > 0$, $I_M = 0$ and there exists no non-zero global minimizers of the free energy.
\end{prop}
\begin{proof}
Recall from above that any global minimizer $u^\star \in \mathcal{Y}_M$ will satisfy $S(u^\star) = \norm{u^\star}_2^2 < \infty$. 
By Cauchy-Schwarz and Young's inequality for convolutions, 
\begin{equation*}
\F(u^\star) = \norm{u^\star}_2^2 - \frac{1}{2}\int\int u^\star(x)u^\star(y) \K(x-y) dxdy \geq \left(1 - \frac{\norm{\K}_1}{2}\right)\norm{u^\star}_2^2. 
\end{equation*}
This is clearly a contradiction unless $\norm{u^\star}_2^2 = 0$, as $I_M \leq 0$ for all $M \geq 0$. 
\end{proof}

We now prove the main theorem. 

\begin{proof}(\textbf{Theorem \ref{thm:Minimizers}})
To prove \textit{(i)} and \textit{(ii)} we begin with the following lemma, which is a restatement of the scaling analysis discussed above. 
\begin{lemma}[Scaling Lemma] \label{lem:infF}
Let \textit{(i)} or \textit{(ii)} hold. Then $\forall \, M > 0$,  $\exists\,\phi \in C^{\infty}_c \cap \mathcal{Y}_M$ with $\F(\phi) < 0$. 
\end{lemma}
\begin{proof}
Let $\phi \in C^{\infty}_c \cap \mathcal{Y}_M$ and consider the mass-invariant scaling $\phi_\lambda(x) = \lambda^d \phi(\lambda x)$. Then for $R > 0$,
\begin{align*}
\F(\phi_\lambda) & = S(\phi_\lambda)-\frac{\lambda^{2d}}{2}\int \int \phi(\lambda x) \phi(\lambda y) \K(x-y) dxdy \\
& = \lambda^d\left(\int \frac{\Phi(\lambda^d\phi(\lambda x))}{\lambda^{d}} dx - \frac{1}{2}\int\int\phi(x)\phi(y)\frac{1}{\lambda^d}\K\left(\frac{x-y}{\lambda}\right) dxdy\right) \\
& \leq \lambda^d\left(\int \frac{\Phi(\lambda^d\phi(\lambda x))}{\lambda^{d}} dx - \frac{1}{2}\int\int\phi(x)\phi(y)\frac{1}{\lambda^d}\K\left(\frac{x-y}{\lambda}\right)\mathbf{1}_{B_R}\left(\abs{x-y}\right) dxdy\right). 
\end{align*}
By the \eqref{def:chi} and $\phi \in C^\infty_c$, for all $\epsilon > 0$ and $\lambda$ sufficiently small (depending on $\phi$) such that,  
\begin{equation*} 
\int \lambda^{-d} \Phi(\lambda^d\phi(\lambda x)) dx = \int \lambda^{-2d} \Phi(\lambda^d\phi(x)) dx \leq \chi \norm{\phi}_2^2 + \epsilon\norm{\phi}_2^2.  
\end{equation*}
Therefore, since $\K \in L_{loc}^1$ and $\phi \in C^\infty_c$, for all $\epsilon$ we may pick $\lambda$ small such that, 
\begin{equation*}
\F(\phi_\lambda) \leq \lambda^d\left(\chi - \frac{\norm{\K\mathbf{1}_{B_R}}_1}{2}\right) \norm{\phi}_2^2 + \epsilon\lambda^d. 
\end{equation*}
Therefore, if $\chi = 0$, we clearly have $I_M < 0$. Moreover,
if $0 < \chi < \infty$, then the sign of the right hand side does not depend on $\lambda$ or $\phi$, and is negative only if
$\norm{\K\mathbf{1}_{B_R}}_1 > 2\chi$. By choosing $R$ sufficiently large, this is equivalent to $\norm{\K}_1 > 2\chi$. 
\end{proof}
Lemma \ref{lem:infF} provides conditions under which $I_M < 0$, but in general this is insufficient to imply the existence of non-trivial minimizers. Without strict sub-additivity, we cannot directly apply the results of \cite{LionsCC84} to attain assertions \textit{(i)} and \textit{(ii)}. To recover, we will use a symmetrization argument and the following lemma, which in general is strictly weaker than sub-additivity. 
\begin{lemma} \label{lem:FE_strict_decrease}
Let (i) or (ii) hold and $M_1 > M_2$. Then $I_{M_1} < I_{M_2}$. 
\end{lemma}
\begin{proof}
Let $u_n \in \mathcal{Y}_{M_2}$ such that $\lim_{n \rightarrow \infty} \F(u_n) \rightarrow I_{M_2}$. By 
Lemma \ref{lem:infF}, there exists $v \in C_c^\infty\cap\mathcal{Y}_M$ such that $\F(v) < 0$ and $\norm{v}_1 = M_1 - M_2$.  
Without loss of generality, by Riesz symmetric decreasing rearrangement and $\K$ radially symmetric non-increasing, we may take $u_n$ and $v$ to be radially symmetric and non-increasing, since applying 
a symmetric rearrangement will only decrease the interaction energy and leaves the entropy unchanged \cite{Lieb83}. Therefore $\forall\,\epsilon > 0$, $\exists\, R_n(\epsilon) > 0$ such that $\norm{u_n\mathbf{1}_{\Real^d\setminus B_{R_n}}}_p < \epsilon$ for $1 \leq p \leq \infty$ (note this does not imply any kind of tightness). Choose $x_n$ such that $B_{R_n}(x_n) \cap \textup{supp} v = \emptyset$ and let $\hat{u}_n = u_n(\cdot - x_n)$. Define,  
\begin{equation} 
z_n(x) = v(x) + \hat{u}_n(x). 
\end{equation}   
By $v,\hat{u}_n \geq 0$ we have $\norm{z_n}_1 = M_1$. Now by the approximately disjoint supports,  
\begin{align*}
\F(z_n) & = \int \Phi(z_n(x)) dx - \frac{1}{2}\int\int z_n(x)z_n(y) \K(x-y) dxdy \\
& \leq \int \Phi(u_n(x)) dx + \int_{\textup{supp}\,v} \Phi(v(x) + \epsilon) dx - \frac{1}{2}\int\int z_n(x)z_n(y) \K(x-y) dxdy.
\end{align*}
Notice that since $v \in C_c^\infty$, by the mean value theorem, 
\begin{equation*}
\int_{\textup{supp}\,v} \Phi(v(x) + \epsilon) dx \leq \int \Phi(v) dx + \mathcal{E}(\epsilon), 
\end{equation*}
such that $\lim\inf_{\epsilon \rightarrow 0}\mathcal{E}(\epsilon) = 0$.  
Therefore, by $\K \geq 0$, 
\begin{equation*}
\F(z_n) \leq \F(u_n) + \F(v) + \mathcal{E}(\epsilon).
\end{equation*} 
Since $\F(v) < 0$ we may choose $\epsilon$ sufficiently small to ensure $\lim\inf_{n \rightarrow \infty}\F(z_n) < \lim\inf_{n \rightarrow \infty}\F(u_n) = I_{M_2}$. 
\end{proof}

We now prove the assertions \textit{(i)} and \textit{(ii)}.
Indeed, let $u_n \in \mathcal{Y}_M$  be such that $\F(u_n) \rightarrow I_M$. 
As above, we may assume that $u_n$ is radially symmetric non-increasing.
We now show $\set{u_n}$ has a convergent subsequence in the strong $L^1$ topology. 

We follow the approach of Theorem II.1 in \cite{LionsCC84}. Following the work contained therein, tightness up to translation is established using the profile decomposition lemma (Lemma I.1). Accordingly, there exists a subsequence of $\set{u_n}$, not relabeled, such that one of the following three possibilities occurs,
\begin{itemize}
\item[(i)] Tight up to translation: $\exists \set{y_n} \subset \Real^d$ for which $\forall\,\epsilon > 0, \; \exists\,R>0$ such that $\int_{\Real^d \setminus B_R(y_n)} u_n dx < \epsilon$. 
\item[(ii)] Vanishing: $\forall R > 0$, $\lim_{n \rightarrow \infty}\sup_{y \in \Real^d} \int_{B_R(y)} u_n dx = 0$. 
\item[(iii)] Dichotomy: $\exists \set{u^1_n},\set{u^2_n},\set{v_n} \subset L_+^1$, such that
  $u_n = u^1_n + u^2_n + v_n$ with, for $i \in \set{1,2}$, $u^1_nu^2_n = u^i_nv_n \equiv 0$, $u_n^i,v_n \leq u_n$,  $\lim_{n \rightarrow \infty} \textup{dist}(\textup{supp}\,u^1_n,\textup{supp}\,u^2_n) = \infty$ and $\lim_{n \rightarrow \infty} \norm{u^1_n}_1 = M - \alpha$, $\lim_{n \rightarrow \infty} \norm{u^2_n}_1 = \alpha$ and $\lim_{n \rightarrow \infty} \norm{v_n}_1 = 0$ for some $\alpha, \; 0 < \alpha < M$. 
\end{itemize}

\textit{Vanishing does not occur:}
\\
Vanishing is ruled out by $I_M < 0$. Indeed, $I_M < 0$ implies $\lim_{n \rightarrow \infty} \mathcal{W}(u_n) > 0$. Assume for contradiction that the subsequence (not relabeled) $\set{u_n}$ vanishes as $n \rightarrow \infty$. Let $q \in [m^\star/(2m^\star - 2), p)$, $R>0$ and by H\"older's inequality,
\begin{align*}
\mathcal{W}(u_n) & = \int\int u_n(x)u_n(y)\K(x-y) dx dy \\
& \leq \norm{u_n}_{2q/(2q-1)}^2\norm{\K\mathbf{1}_{B_{R^{-1}}}}_q + \int\int_{R^{-1} < \abs{x-y} \leq R}u_n(x)u_n(y) \K(x-y) dx dy + M^2\norm{\K\mathbf{1}_{\Real^d\setminus B_R}}_\infty \\
& \leq \norm{u_n}_{2q/(2q-1)}^2\norm{\K\mathbf{1}_{B_{R^{-1}}}}_q + \norm{\K \mathbf{1}_{B_{R}\setminus B_{R^{-1}}}}_\infty\int u_n(x) \int_{\abs{x-y}<R} u_n(y) dy dx + M^2\norm{\K\mathbf{1}_{\Real^d\setminus B_R}}_\infty.
\end{align*}
By interpolation, $\norm{u_n}_{2q/(2q-1)}$ is uniformly bounded by $2q/(2q-1) \leq m^\star$, and since $\set{u_n}$ vanishes we may deduce, 
\begin{align*}
\lim\inf_{n \rightarrow \infty} \mathcal{W}(u_n) \lesssim \norm{\K\mathbf{1}_{B_{R^{-1}}}}_q + \norm{\K\mathbf{1}_{\Real^d\setminus B_R}}_\infty.  
\end{align*}
As $R \rightarrow \infty$, the last term vanishes since $\K$ is radially symmetric non-increasing and the first 
vanishes by the dominated convergence theorem and $\K \in L_{loc}^q$. Therefore, we have deduced
\begin{equation*}
\lim\inf_{n \rightarrow \infty} \mathcal{W}(u_n) \leq 0, 
\end{equation*}
which is a contradiction to $I_M < 0$. 

\textit{Dichotomy does not occur:}
\\
Although we do not have strict sub-additivity, we will take advantage of the weaker property, Lemma \ref{lem:FE_strict_decrease}, along with radial symmetry, to rule out dichotomy. Suppose for contradiction that dichotomy occurs. 
By Riesz symmetric decreasing rearrangement, recall that $u_n$ is radially symmetric non-increasing. This together with the properties of the profile decomposition $u_n = u_n^1 + u_n^2 + v_n$ implies
one of $u^1_n$ or $u^2_n$ converges to zero in $L^\infty$. In particular, one of the sequences must vanish, which is the advantage of radial symmetry.
Assume without loss of generality that $u^2_n \rightarrow 0$. First note by the disjoint supports, 
\begin{equation}
S(u_n) \geq S(u_n^1) + S(u_n^2). \label{ineq:entropy}
\end{equation}
Then, 
\begin{align*}
\mathcal{W}(u_n) & = \mathcal{W}(u_n^1) + \mathcal{W}(u_n^2) + \mathcal{W}(v_n) \\ 
& + \int\int u_n^1u_n^2\K(x-y)dxdy + \int\int v_nu_n^1\K(x-y)dxdy + \int\int v_nu_n^2\K(x-y)dxdy.
\end{align*}
By the \eqref{ineq:GHLS} and interpolation, for any $\delta > 0$, 
\begin{align}
\mathcal{W}(u_n^2) & \lesssim \wknorm{\K\mathbf{1}_{B_\delta}}{p}\norm{u_n^2}_1^{2-m^\star}\norm{u_n^2}_{m^\star}^{m^\star} + \norm{\K\mathbf{1}_{\Real^d\setminus B_\delta}}_{\hat{p}}\norm{u_n^2}^2_{\frac{2\hat{p}}{2\hat{p} - 1}} \nonumber \\ 
& \lesssim  \wknorm{\K\mathbf{1}_{B_\delta}}{p}\norm{u_n^2}_{\infty}^{m^\star - 1} + \norm{\K\mathbf{1}_{\Real^d\setminus B_\delta}}_{\hat{p}}\norm{u_n^2}^{\frac{1}{\hat{p}}}_\infty. \label{ineq:u2}
\end{align}
Similarly for any $\delta > 0$,
\begin{equation}
\mathcal{W}(v_n) \leq \norm{v_n}_1^2\norm{\K\mathbf{1}_{\Real^d \setminus B_\delta}}_\infty + \wknorm{\K\mathbf{1}_{B_\delta}}{p}\norm{v_n}_1^{2-m^\star}\norm{v_n}_{m^\star}^{m^\star}, \label{ineq:vv}
\end{equation}
and for $i \in \set{1,2}$,
\begin{equation}
\int\int v_nu_n^i\K(x-y)dxdy \leq \norm{v_n}_1\norm{u_n^i}_1\norm{\K\mathbf{1}_{\Real^d \setminus B_\delta}}_\infty + \wknorm{\K\mathbf{1}_{B_\delta}}{p}\norm{v_n}_1^{1-m^\star/2}\norm{v_n}_{m^\star}^{m^\star/2}\norm{u_n^i}_1^{1-m^\star/2}\norm{u_n^i}_{m^\star}^{m^\star/2}.\label{ineq:vui}
\end{equation}
Finally, let $d_n = \textup{dist}(\textup{supp}\,u^1_n,\textup{supp}\,u^2_n)$. Therefore, 
\begin{align}
\int\int u_n^1(x)u_n^2(y)\K(x-y)dxdy & = \int\int u_n^1(x)u_n^2(y) \K(x-y) \mathbf{1}_{\Real^d\setminus B_{d_n}}(\abs{x-y}) dx dy \nonumber \\
 & \leq M^2\norm{\K \mathbf{1}_{\Real^d \setminus B_{d_n}}}_\infty. \label{ineq:u1u2}
\end{align}
Putting the estimates \eqref{ineq:entropy}-\eqref{ineq:u1u2} together with $d_n \rightarrow \infty$, 
$\K$ radially symmetric non-increasing, along with $\lim_{n \rightarrow \infty} \norm{u_n^2}_\infty =\lim_{n \rightarrow \infty} S(u_n^2) = \lim_{n \rightarrow \infty} \norm{v_n}_1 = 0$ and the uniform boundedness of $\norm{u_n^i}_{m^\star},\norm{v_n}_{m^\star}$ implies, 
\begin{align*}
I_M = \lim_{n \rightarrow \infty} \F(u_n) \geq \lim\inf_{n \rightarrow \infty}\F(u_n^1) \geq I_{M-\alpha}.
\end{align*}
This clearly contradicts Lemma \ref{lem:FE_strict_decrease} and rules out dichotomy, leaving 
only that there is a subsequence of $\set{u_n}$ which is tight up to translation.  

\textit{Conclusion of proof:} \\
Following \cite{LionsCC84} one may now prove without modification that tightness is sufficient 
to extract a subsequence (not relabeled) such that $u_n \rightarrow u^\star$ strongly in $L^1$ for some $u^\star \in \mathcal{Y}_M$ with $\F(u^\star) = I_M$. Strong convergence is due to $\F(u_n) \rightarrow I_M$ and the strict convexity of $\Phi(u)$ (see \cite{LionsCC84}). This concludes the proof of assertions \textit{(i)} and \textit{(ii)}. 
\end{proof}

\section{Acknowledgments} 
The author would like to thank Andrea Bertozzi, Inwon Kim, Thomas Laurent, Nancy Rodr\'iguez and Yao Yao for their guidance and helpful discussions.
This work was in part supported by NSF grant DMS-0907931.

\bibliographystyle{plain}
\bibliography{nonlocal_eqns}

\begin{thebibliography}{10}

\bibitem{AmbrosioGigliSavare}
L.A. Ambrosio, N.~Gigli, and G.~Sava\'re.
\newblock {\em Gradient flows in metric spaces and in the space of probability
  measures}.
\newblock Lectures in Mathematics, Birkh\"auser, 2005.

\bibitem{BRB10}
J.~Bedrossian, N.~Rodr\'iguez, and A.L. Bertozzi.
\newblock Local and global well-posedness for aggregation equations and
  {Patlak-Keller-Segel} models with degenerate diffusion.
\newblock {\em Preprint, {\textup{arXiv:1009.2674}}}, 2010.

\bibitem{BertozziSlepcev10}
A.L. Bertozzi and D.~Slepcev.
\newblock Existence and uniqueness of solutions to an aggregation equation with
  degenerate diffusion.
\newblock {\em To appear in Comm. Pure. Appl. Anal.}, 2010.

\bibitem{BlanchetCalvezCarrillo08}
A.~Blanchet, V.~Calvez, and J.A. Carrillo.
\newblock Convergence of the mass-transport steepest descent scheme for
  subcritical {Patlak-Keller-Segel} model.
\newblock {\em SIAM J. Num. Anal.}, 46:691--721, 2008.

\bibitem{Blanchet09}
A.~Blanchet, J.~Carrillo, and P.~Laurencot.
\newblock Critical mass for a {Patlak-Keller-Segel} model with degenerate
  diffusion in higher dimensions.
\newblock {\em Calc. Var.}, 35:133--168, 2009.

\bibitem{Burger07}
M.~Burger, V.~Capasso, and D.~Morale.
\newblock On an aggregation model with long and short range interactions.
\newblock {\em Nonlin. Anal. Real World Appl.}, 8(3):939--958, 2007.

\bibitem{CarrilloEntDiss01}
J.A. Carrillo, A.~J\"ungel, P.A. Markowich, G.~Toscani, and A.~Unterreiter.
\newblock Entropy dissipation methods for degenerate parabolic problems and
  generalized {Sobolev} inequalities.
\newblock {\em Montash. Math.}, 133:1--82, 2001.

\bibitem{FellnerRaoul1}
K.~Fellner and G.~Raoul.
\newblock Stability of stationary states of non-local equations with singular
  interaction potentials.
\newblock {\em to appear in {Math. Comp. Model.}}, 2010.

\bibitem{FellnerRaoul2}
K.~Fellner and G.~Raoul.
\newblock Stable stationary states of non-local interaction equations.
\newblock {\em to appear in {Math. Mod. Meth. Appl. Sci.}}, 2010.

\bibitem{YaoKim10}
I.~Kim and Y.~Yao.
\newblock {\em personal communication}.

\bibitem{LiRodrigo09}
D.~Li and J.~Rodrigo.
\newblock Finite-time singularities of an aggregation equation in
  {$\mathbb{R}^n$} with fractional dissipation.
\newblock {\em Comm. Math. Phys.}, 287:687--703, 2009.

\bibitem{Lieb83}
E.H. Lieb.
\newblock Sharp constants in the {Hardy-Littlewood-Sobolev} and related
  inequalities.
\newblock {\em Ann. Math.}, 118:349--374, 1983.

\bibitem{LiebYau87}
E.H. Lieb and H-T Yau.
\newblock The {Chandrasekhar} theory of stellar collapse a the limit of quantum
  mechanics.
\newblock {\em Comm. Math. Phys.}, 112:147--174, 1987.

\bibitem{LionsCC84}
P.L. Lions.
\newblock The concentration-compactness principle in calculus of variations.
  the locally compact case, part 1.
\newblock {\em Ann. Inst. Henri. Poincare, Anal. non lin}, 1(2):109--145, 1984.

\bibitem{McCann97}
R.J. McCann.
\newblock A convexity principle for interacting gases.
\newblock {\em Adv. Math.}, 128:153--179, 1997.

\bibitem{Otto01}
F.~Otto.
\newblock The geometry of dissipative evolution equations: the porous medium
  equation.
\newblock {\em Comm. Part. Diff. Eqn.}, 26(1):101--174, 2001.

\bibitem{Slepcev08}
D.~Slep\v{c}ev.
\newblock Coarsening in nonlocal interfacial systems.
\newblock {\em SIAM J. Math. Anal.}, 40(3):1029--1048, 2008.

\bibitem{SugiyamaNA05}
Y.~Sugiyama.
\newblock Global existence and decay properties of solutions for some
  degenerate quasilinear parabolic systems modelling chemotaxis.
\newblock {\em Nonlinearity}, 63:1051--1062, 2005.

\bibitem{SugiyamaDIE06}
Y.~Sugiyama.
\newblock Global existence in sub-critical cases and finite time blow-up in
  super-critical cases to degenerate {Keller-Segel} systems.
\newblock {\em Diff. Int. Eqns.}, 19(8):841--876, 2006.

\bibitem{SugiyamaADE07}
Y.~Sugiyama.
\newblock Application of the best constant of the {Sobolev} inequality to
  degenerate {Keller-Segel} models.
\newblock {\em Adv. Diff. Eqns.}, 12(2):121--144, 2007.

\bibitem{SugiyamaDIE07}
Y.~Sugiyama.
\newblock The global existence and asymptotic behavior of solutions to
  degenerate to quasi-linear parabolic systems of chemotaxis.
\newblock {\em Diff. Int. Eqns.}, 20(2):133--180, 2007.

\bibitem{TopazBertozzi04}
C.M. Topaz and A.L. Bertozzi.
\newblock Swarming patterns in a two-dimensional kinematic model for biological
  groups.
\newblock {\em SIAM J. Appl. Math.}, 65(1):152--174, 2004.

\bibitem{TopazBertozziLewis06}
C.M. Topaz, A.L. Bertozzi, and M.A. Lewis.
\newblock A nonlocal continuum model for biological aggregation.
\newblock {\em Bull. Math. Bio}, 68(68):1601--1623, 2006.

\end{thebibliography}

\end{document}